\documentclass{amsproc}

\usepackage{a4wide}
\usepackage{amsmath}
\usepackage{enumerate}

\usepackage{amsmath,amsthm}
\usepackage{amsfonts}
\usepackage{amssymb}
\usepackage{longtable}
\usepackage[matrix,arrow,curve]{xy}

\newtheorem{theorem}{Theorem}
\newtheorem{cor}{Corollary}

\newtheorem{lem}{Lemma}

{\bf}{\rm}

\newtheorem{rem}{Remark}{\bf}{\rm}{\rm}

\newtheorem*{theorem*}{Theorem}
\newtheorem*{cor*}{Corollary}

\def\Real{\mathbb{R}}

\def\D{\mathcal{D}}
\def\C{\mathcal{C}}

\def\F{\mathcal{F}}
\def\vH{{\check{H}}}

\def\On{\text{\rm O}_n}

\def\GL1{\text{\rm GL}(1,\Real)}
\def\O1{\text{\rm O}_1}

\def\GVL{\text{\rm GVL}}

\def\id{\mathop\text{\rm id}\nolimits}
\def\pt{\mathop\text{\rm pt}\nolimits}

\usepackage{stackrel}

\newcommand{\be}{\begin{equation}}
\newcommand{\ee}{\end{equation}}

\let\ge=\geqslant
\let\leq=\leqslant
\let\geq=\geqslant

\begin{document}

\title{Non-diffeomorphic Reeb foliations and modified Godbillon-Vey class}

\author{Yaroslav V. Bazaikin}\thanks{$^1$Sobolev Institute of Mathematics, Novosibirsk, Russia and University of Hradec Kr\'alov\'e, Faculty of Science, Rokitansk\'eho 62, 500~03 Hradec Kr\'alov\'e,  Czech Republic
E-mail: bazaikin(at)math.nsc.ru}

\author{Anton S. Galaev}\thanks{$^2$University of Hradec Kr\'alov\'e, Faculty of Science, Rokitansk\'eho 62, 500~03 Hradec Kr\'alov\'e,  Czech
Republic\\
E-mail: anton.galaev(at)uhk.cz}

\author{Pavel Gumenyuk}\thanks{$^3$ Department of Mathematics, Politecnico di Milano, via E.~Bonardi 9, 20133 Milan, Italy.}

\maketitle

\begin{abstract}

The paper deals with a modified Godbillon-Vey class defined by Losik for codimension-one foliations. This characteristic class takes values in the cohomology of the second order frame bundle over the leaf space of the foliation. 
The definition of the Reeb foliation depends upon two real functions satisfying certain conditions. All these foliations are pairwise homeomorphic and have trivial Godbillon-Vey class. We show that the 
modified Godbillon-Vey  is non-trivial for some Reeb foliations and it is trivial for some other Reeb foliations. In particular, the modified Godbillon-Vey class can distinguish non-diffeomorphic foliations and it provides more information than the classical Godbillon-Vey class. 
 We also show that this class is non-trivial for some foliations on the two-dimensional surfaces.

% We construct explicit examples of the Reeb foliations that are not diffeomorphic. 

\vskip0.5cm

{\bf Keywords}: Reeb foliation; Reeb component; leaf space of foliation; characteristic classes of foliation; Gelfand formal geometry; Gelfand-Fuchs cohomology; Godbillon-Vey-Losik
class.

\vskip0.5cm

{\bf AMS Mathematics Subject Classification:} 57R30; 57R32.

\tableofcontents

%    57R30   Foliations; geometric theory
%57R32   Classifying spaces for foliations; Gelfand-Fuks cohomology

\end{abstract}

\section*{Introduction}

The present paper is a continuation of our recent work \cite{BG} on characteristic classes of codimension-one foliations based on Losik's ideas~\cite{L15,L90,L94} that appeared as the development of  Gelfand's formal geometry, cf. \cite{Kontsevich}. For a codimension-one foliation $\F$ on a smooth manifold $M$, the following sequence of homomorphisms is defined \cite{BG,CM,Gal17}:
$$H^3(W_1,\O1)\to H^3(S_2(M/\F)/(\O1\times\mathbb{Z}))\to H^3(S_2(M/\F)/\O1)\to \vH^3(M/\F)\to H^3(M),$$ where $H^3(W_1,\O1)$  is the Gelfand-Fuchs cohomology of the Lie algebra of formal vector fields on the line, $S_2(M/\F)$ is the second order frame bundle over the leaf space 
 $M/\F$, and $\vH^3(M/\F)$ is the \v{C}ech-de\,Rham cohomology of the leaf space studied in details in \cite{CM}. The generator of $H^3(W_1,\O1)\cong\Real$ defines the classical Godbillon-Vey class with values in $H^3(M)$ as well as three modifications of that class. 

The classical Godbillon-Vey class is an important invariant of a codimension-one foliation, see, e.g., the survey \cite{Hurder00} and the references therein.   The Godbillon-Vey number is an invariant for the cobordism classes of foliations on three dimensional manifolds and in some cases it allows to distinguish non-cobordant and, consequently, non-diffeomorphic foliations \cite{Ref106}. In \cite{CM} is given an example of a foliation such that its   Godbillon-Vey class in non-trivial in $\vH^3(M/\F)$ and it is trivial in $H^3(M)$.

We fix our attention on the Reeb foliations, which play a significant role in the theory of codimension-one foliations. They are the elementary ``building blocks'' filling the gaps between other foliation components. For instance, from the  celebrated Novikov compact leaf theorem it follows that each foliation on  $S^3$  contains a Reeb component. According to the classification theory of foliations on  two-dimensional surfaces, the  Reeb component may be always embedded in a neighbourhood of a closed leaf.

The classical Godbillon-Vey class is trivial for all Reeb foliations. 
This may be checked directly~\cite{Tamura}. This  also  follows  from~\cite{Ref82} since the Reeb foliations are almost without holonomy. The Godbillon-Vey class with values in $\vH^3(M)$ is also trivial for all Reeb foliations, this follows from the results in Noncommutative geometry \cite{Connes}. In \cite{BG} (see also \cite{L15}), we show that the Godbillon-Vey class with value in  
$H^3(S_2(M/\F)/(\O1\times\mathbb{Z}))$ is non-trivial for all Reeb foliations. That, in particular, implies that this class detects the compact leaf with non-trivial holonomy; this is the information that cannot be obtained neither from  the Godbillon-Vey class with values in $\vH^3(M/\F)$ nor from the classical Godbillon-Vey class. 

In the present paper we fix our attention on the Godbillon-Vey class with value in $H^3(S_2(M/\F)/\O1)$. We call this class the Godbillon-Vey-Losik class (GVL class). It turns out that this class is very sensitive to the dynamics of the non-compact leaves in the following sense. Dynamics of the non-compact leaves in a neighborhood of the compact leaf $\mathcal {L}$ is described by the holonomy group of  $\mathcal {L}$. For a fixed $x\in\mathcal{L}$ and a transversal $T$ to $\F$ through $x$, the holonomy group of $\mathcal{L}$  consists of the germs at $x$ of local diffeomorphisms  of  $T$  defined by loops in $\mathcal {L}$ starting at $x$. This group can be included to a $1$-parameter group of the germs of local diffeomorphisms of the transversal with one fixed point $x$. Even if this group consists of the germs of diffeomorphisms infinitely tangent to the identity, there is a notion of the order of the convergence (a $1$-parameter group has greater order of convergence than another $1$-parameter group, if the fraction of their generating vector fields is smooth and has zero value at $x$). It turns out that in some situations the GVL class distinguishes this order. Now we state the main results of the paper.

To define the Reeb foliation on the solid torus, one fixes an even function $f:(-1,1)\to\Real$ with certain properties. In the present paper we consider functions $f$ that are roughly speaking given by the condition $f'(x)=e^\frac{1}{(1-x)^\alpha}$ in a neighbourhood of $x=1$, here $\alpha\in \Real$, $\alpha>0$. Let $\mathcal{R}^+_\alpha$ denote the corresponding Reeb foliation on the solid torus.

\begin{theorem*}
	The GVL class of the foliation ${\mathcal R}^+_\alpha$ is trivial if and only if $\alpha \in \mathbb{N}$.
\end{theorem*}

\begin{cor*} If $\alpha \in \mathbb{N}$ and $\beta \notin \mathbb{N}$, then the
	 Reeb foliations ${\mathcal R}_\alpha^+$ and ${\mathcal R}_\beta^+$ are not diffeomorphic. 	
\end{cor*}

One can glue two Reeb foliations $\mathcal{R}^+_\alpha$  to obtain a Reeb foliation ${\mathcal R}_\alpha$ on the sphere $S^3$.

\begin{theorem*}
	The GVL class of ${\mathcal R}_\alpha$ is trivial if and only if $\alpha \in 2\mathbb{N}$.
\end{theorem*}

\begin{cor*} If $\alpha \in 2\mathbb{N}$ and $\beta \notin 2\mathbb{N}$, then
the	Reeb foliations ${\mathcal R}_\alpha$ and ${\mathcal R}_\beta$ are not diffeomorphic.
\end{cor*}

 Mizutani \cite{Mizutani} and Sergertaert \cite{Serg} proved that any Reeb foliation is cobordant to zero. It follows that ${\mathcal R}_\alpha$ and ${\mathcal R}_\beta$ are cobordant and in contrast to the properties of the classical Godbillon-Vey classes we have:

\begin{cor*}
The	GVL class is not a cobordism invariant.
\end{cor*}

We consider also the GVL class for the Reeb foliations on the ring $[-1,1]\times S^1$ and obtain for it similar results as for the foliations on the solid torus.

\section{Definition of characteristic classes following Losik}\label{sec1}

Here we review Losik's approach to the leaf spaces of foliations
\cite{L90,L94,L15} that allowed him to define new characteristic
classes of foliations. Let $\D_n$ be the category whose objects
are open subsets of $\Real^n$, and morphisms are \'etale (i.e.,
regular) maps.

Let us recall the definition of a $\D_n$-space. Let $X$ be a set. A
$\D_n$-chart on $X$ is a pair $(U,k)$, where $U\subset \Real^n$ is an open subset, and
$k:U\to X$ is an arbitrary map. For two charts $k_i:U_i\to X$, a
morphism of charts is an \'etale map $m:U_1\to U_2$ such that
$k_2\circ m=k_1$. Let $\Phi$ be a set of charts and let $\C_\Phi$ be
a category whose objects are elements of $\Phi$ and morphisms are
some morphisms of the charts. The set $\Phi$ is called a
$\D_n$-atlas on $X$ if $X=\varinjlim J$, where $J:C_\Phi\to{\rm
Sets}$ is the obvious functor. A $\D_n$-space is a set $X$ with a maximal $\D_n$-atlas $\Phi$. A morphism of $\D_n$ spaces is a map $\varphi:X\to Y$ such that for each chart $k:U\to X$ from the $\D_n$-atlas on $X$, the map $\varphi\circ k:U\to Y$ is a chart from the maximal $\D_n$-atlas on $Y$.

Let $X$ be a $\D_n$-space. To define a $q$-form on $X$,  one fixes a $q$-form $\omega_U$ on $U$ for each chart $k:U\to X$ in such a way that for each morphism $m:U_1\to U_2$ it holds $m^*\omega_{U_2}=\omega_{U_1}.$ The exterior derivative is defined on such forms and we obtain the de~Rham cohomology $H^*(X)$.

If $\F$ is a foliation of codimension $n$ on a smooth manifold $M$,
then the leaf space $M/\F$ is a $\D_n$-space. The maximal
$\D_n$-atlas on $M/\F$ consists of the projections $U\to M/\F$,
where $U$ is a transversal which is an open subset
of $\Real^n$ embedded to $M$. These transversals may be obtained from a foliation
atlas on $M$.

Generally $\D_n$-spaces are orbit spaces of  pseudogroups of local
diffeomorphisms of smooth manifolds. Considering the space $\Real^n$
and the pseudogroup of all local diffeomorphisms of open subsets of
$\Real^n$, we see that the point $\pt$ is a $\D_n$-space. The atlas
of $\pt$ consists of all pairs $(U,k)$, where $U\subset\Real^n$ is an
open subset and $k:U\to \pt$ is the unique map. It is important to
note that $\pt$ is the terminal objects in the category of
$\D_n$-spaces.

Let $X$ be a $\D_n$ space. For each open subset $U\subset\Real^n$ consider the space $S_2(U)$ of the second order jets   of regular maps $f:W\to U$ at zero, where $W\subset \Real^n$ is an open subset containing the origin. Consider also the factor $=S_2(U)/\On.$
If $k_i:U_i\to X$, $i=1,2$, are $\D_n$-charts and $m:U_1\to U_2$ is a morphism of charts, then we get the induced morphism $$\tilde m:S_2(U_1)/\On\to S_2(U_2)/\On.$$ The sets $S_2(U)/\On$ and morphisms $\tilde m$ define a $\D_{N}$-space (for some number $N$), which we denote by $S_2(X)/\On$.
For the terminal object $\pt$ it holds  $$H^*(S_2(\pt)/\On)\cong H^*(W_n,\On),$$ where
$H^*(W_n,\On)$ are the  relative Gelfand-Fuchs cohomology of the Lie algebra of formal vector fields on $\Real^n$. The unique map $X\to\pt$ induces the characteristic homomorphism \begin{equation}
\label{charmap}H^*(W_n,\On)\to  H^*(S_2(X)/\On).\end{equation}
The image of the generators of $H^*(W_n,\On)$ under that map are characteristic classes of the $\D_n$-space~$X$. If $X=M/\F$ is the leaf space of a
 codimension $n$ foliation, then there is a homomorphism
\begin{equation}\label{prtoM} H^*(S_2(X)/\On)\to H^*(M), \end{equation}
which maps these classes to the usual characteristic classes of the foliation $\F$. We will see that this homomorphism may be non-injective.

Let now $n=1$. Let $U\subset \Real$ be an open subset.  Consider the coordinates
$$
 z_p =\frac{d^p f}{(dt)^p}(0),\quad p=0,1,2$$
 on $S_2(U)$, where $f:W\to U$ is a regular smooth map, and $W\subset \Real$ is open and contains the origin. We will consider the following coordinates
 on $S_2(U)/\On$:
\begin{equation} \label{coord3} x_0 =z_0,\quad x_1=\ln|z_1|,\quad x_2=\frac{z_2}{z_1^2}. \end{equation}
If $f:U\to V$, $y=f(x)$, is a local diffeomorphis, then the map $\tilde f:S''_2(U)\to S''_2(V)$ is given by
\begin{equation}\label{coordchange}y_0=f(x_0),\quad y_1=x_1+\ln|f'(x_0)|,\quad
y_2=\frac{x_2}{f'(x_0)}+\frac{f''(x_0)}{(f'(x_0))^2}.\end{equation}
The form $$-dx_0\wedge dx_1\wedge d x_2$$
is invariant with respect to any such coordinate transformation, i.e. it defines 3-form on $S_2(X)/\On$. The cohomology class of this form is the image of the canonical generator of $H^3(W_1,\O1)\cong\Real$ under the map \eqref{charmap}. We call this class the {\bf Godbillon-Vey-Losik class} (GVL class). If $X=M/\F$, then the map \eqref{prtoM} maps this class to the classical Godbillon-Vey class.

The GVL class of a $\D_1$-space is zero if and only if for each chart $k:U\to X$ a 2-from $\omega_U$ on $U$ with $$d\omega=-dx_0\wedge dx_1\wedge d x_2$$ is fixed, and for each morphism $m:U_1\to U_2$ it holds $$\tilde m^*\omega_{U_2}=\omega_{U_1}.$$

The definition of the GVL class is extended naturally to codimension-one foliations on manifolds with boundaries (the foliations are tangent to the boundary). In that case, the domains of the charts may be semi-closed intervals, and we assume that a function is smooth on a subset of $\Real^N$ if it may be smoothly extended to an open neighborhood of that subset.

\section[]{Reeb foliations}\label{secReeb}

The definition of the Reeb foliations may be found in many texts, see, e.g., \cite{Tamura}.
Consider the three-dimensional sphere $S^3$ of radius $\sqrt{2}$,
\begin{equation*}
S^3 = \{ (z_1, z_2) : z_1, z_2, \in \mathbb{C}, |z_1|^2 + |z_2|^2 = 2 \}.	
\end{equation*}
Let $$g, f: (-1, 1) \rightarrow \mathbb{R}$$ be two smooth even functions strictly increasing on $(0,1)$ and satisfying the conditions \begin{align} \lim_{t\to\pm 1}g(t)&=\lim_{t\to\pm 1}f(t)=+\infty,\\ \label{propf2}
\lim_{t\to\pm 1}\left(\frac{1}{g'(t)}\right)^{(k)}&=\lim_{t\to\pm 1}\left(\frac{1}{f'(t)}\right)^{(k)}=0, \quad k=0,1,2,\dots.\end{align}
Note that in some texts, the condition $\lim_{t\to\pm 1}g^{(k)}(t)=\lim_{t\to\pm 1}f^{(k)}(t)=+\infty$, $k=1,2,\dots$, is also assumed.
The
Reeb foliation ${\mathcal R}_{g,f}$ on $S^3$ is defined by the following leaves:
\begin{align*}
L_\xi^- &= \left\{ \left(z, \sqrt{2-|z|^2}e^{2\pi i(g(|z|) - \xi)}\right)  : z \in \mathbb{C}, |z|<1  \right\}, \\
L_\zeta^+ &= \left\{ \left(\sqrt{2-|z|^2}e^{2\pi i(f(|z|)- \zeta)}, z \right) :  z \in \mathbb{C}, |z|<1  \right\}, \\
L_0 &= \left\{ \left( e^{2\pi i a}, e^{2\pi i b} \right) : a,b \in \mathbb{R} \right\},
\end{align*}
where $\xi, \zeta \in \mathbb{R}$ are periodic parameters: $L_\xi^- = L_{\xi + 1}^-$, $L_\zeta^+ = L_{\zeta + 1}^+$. The leaves $L_\xi^-$ and $L_\zeta^+$ are diffeomorphic to the two-dimensional disk $D^2 = \{z:\, |z|<1\} \subset \mathbb{C}$ and converge to the singular leaf $L_0$, which is diffeomorphic to the two-dimensional torus.

Define the transversal $$h: \Real \rightarrow S^3$$ of the foliation ${\mathcal R}_{g,f}$ as it follows:
\begin{equation*}
h(x) = \left( \sqrt{2} \sin \left( \frac{\pi}{4} +\frac{1}{2}\arctan x\right), \sqrt{2} \cos \left( \frac{\pi}{4} +\frac{1}{2}\arctan x\right) \right) \in S^3, \quad x \in \mathbb{R}.	
\end{equation*}
Then
\begin{align*}
h(x) &\in L^-_\xi \text{ for } \xi =\hat{g}(x) = g\left(\sqrt{2} \sin \left(\frac{\pi}{4} +\frac{1}{2}\arctan x\right)\right),& x<0, &&   \\
h(x) &\in L^+_\zeta \text{ for } \zeta = \hat{f}(x) = f\left(\sqrt{2} \cos \left(\frac{\pi}{4} +\frac{1}{2}\arctan x\right)\right),& x>0, && \\
h(0) &\in L_0.	
\end{align*}

The fundamental group of the torus $L_0$ is generated by the following two loops: $$\gamma_1(t) = (1, e^{2\pi i t}),\quad \gamma_2(t) = (e^{2\pi i t},1), \quad 0\leq t \leq 1.$$ Then the holonomy of the singular leaf $L_0$ is generated by the following  local diffeomorphisms $\varphi$ and $\psi$ of the transversal $\Real$:

\begin{equation*}
\begin{aligned}
\varphi(x) = \left\{
\begin{array}{lr}
\hat{f}^{-1}(\hat{f}(x) + 1), & x>0, \\
x, & x \leq 0,
\end{array} \right.
\\
\psi(x) = \left\{
\begin{array}{lr}
\hat{g}^{-1}(\hat{g}(x) + 1), & x < 0, \\
x, & x \geq 0.
\end{array} \right.
\end{aligned}
\end{equation*}

 There are two components of the Reeb foliation which are the restrictions of ${\mathcal R}_{g,f}$ to the two solid tori $T^- = \{(z_1,z_2): |z_1| \leq 1 \} \subset S^3$ with the leaves $L_\xi^-$, $L_0$, and $T^+ = \{(z_1,z_2): |z_2| \leq 1 \} \subset S^3$ with the leaves $L_\zeta^+$, $L_0$. We denote them by ${\mathcal R}^-_g$ and~${\mathcal R}^+_f$, respectively.

We see now that the leaf space $X=S^3/{\mathcal R}_{g,f}$ of the Reeb foliation ${\mathcal R}_{g,f}$ is the orbit space $$\Real/<\psi,\varphi>$$ of the pseudogroup generated by the local diffeomorphisms  $\psi$ and $\varphi$ of $\Real$.  The GVL class of the $\D_1$-space $S^3/{\mathcal R}_{g,f}$ is trivial if and only if there exists a 2-form $\omega$ on $\Real^3=S_2\big(\Real\big)/\O1$ invariant with respect to $\tilde\psi$ and $\tilde\varphi$, and such that with respect to the coordinates \eqref{coord3} it holds
$$d\omega=-dx_0\wedge dx_1\wedge d x_2.$$
Likewise, in the case of a foliation  ${\mathcal R}^+_f$ on the solid torus, the form $\omega$ must be smooth on  $$\mathbb{R}^3_+ = \{ (x_0, x_1, x_2)\, :\, x_i \in \mathbb{R}, x_0\geq 0\}$$ and  $\tilde\varphi$-invariant.

\begin{rem}\label{rem} It is important to note that the leaf space $X=S^3/{\mathcal R}_{g,f}$ may be considered also as the orbit space $(-\varepsilon,\varepsilon)/<\psi,\varphi>$ of the pseudogroup defined by the restrictions of  $\psi$ and $\varphi$, here  $\varepsilon$ is an arbitrary positive real number. By that reason, the triviality of the GVL class is equivalent to the existence of the form $\omega$ defined for $x_0$ from a neighborhood of $0$ in $\Real$.
\end{rem}

Thanks to the above remark, in what follows we may suppose that $\hat f$ and $\hat g$ are modified outside some neighborhood of~$x=0$ in such a way that $\varphi$ and $\psi$ become (global) diffeomorphisms of~$\Real$.

Let $\varphi$ be any  diffeomorphism of $\Real$ such that $\varphi(x)=0$ for $x\leq 0$.
From a general theory it is known  \cite{Kopell,Szekeres} that $\varphi$ can be extended to a one-parameter group of $C^2$-diffeomorphisms $\varphi_t$, $t\in\Real$, with $\varphi_1 = \varphi$,  generated by a vector field $V$ of class~$C^1$ . The vector field $V$ is called the Szekeres vector field. In our case, this one-parameter group can be constructed explicitly:
$$
\varphi_t(x) = \left\{
\begin{array}{lr}
\hat{f}^{-1}(\hat{f}(x) + t), & x>0, \\
x, & x \leq 0.
\end{array} \right.
$$
For $x>0$, it holds $V(x)=1/\hat{f}'(x)$. Conditions \eqref{propf2} imply that
$V$ can be extended smoothly to~$\Real$,
\begin{equation}
V(x) = \left\{
\begin{array}{cc}
1/\hat{f}'(x), & x>0, \\
0, & x \leq 0,
\end{array}\right.
\end{equation}
and, consequently, the one-parameter group $\varphi_t$ is smooth.

Conversely, if we consider any $C^\infty$-function $V(x)$ such that $V(x)=0$ for $x\leq 0$, and $V(x)<0$ for $x>0$, then we can find the function $f(t)$ for $t$ from an interval $(1-\varepsilon,1)$:
$$f(t) = \hat{f}\left(\tan \left( 2\arccos \frac{t}{\sqrt{2}} - \frac{\pi}{2}\right)\right),\quad \hat{f} (x) = \int_1^x \frac{d\xi}{V(\xi)},$$
and extend $f$ to a smooth even function on $(-1,1)$ strictly increasing on $(0,1)$. The geometry of the Reeb foliation depends only on the behavior of the function $f$ in a neighborhood of $t=1$, and, in particular, the leaf space does not depend on the extension of $f$ to $(-1,1)$.

In a similar way one can define the Szekeres field~$\tilde V$ associated to the  diffeomorphism $\psi$ with the support contained in the negative real numbers. This field allows to recover the function $g$.

\section[]{Reeb foliations with non-trivial GVL class}

For real $\alpha>0$ consider the following smooth Szekeres vector field on $\mathbb{R}$:
\begin{equation}\label{Eq4.1}
V(x) =
\left\{ \begin{array}{ll}
-e^{-\frac{1}{x^\alpha}}, & \mbox{for } x >0, \\
0, & \mbox{for }  x\leq 0. 	
\end{array}\right.
\end{equation}
Let $f$ be the function defined by the field $V$. We denote the Reeb component ${\mathcal R}^+_f$ by ${\mathcal R}^+_\alpha$.

\begin{theorem}\label{Th4.1}
	If $\alpha \notin \mathbb{N}$, then the class $\GVL({\mathcal R}^+_\alpha)$ is non-trivial.
\end{theorem}

{\bf Proof}. Let for a moment $V$ be the Szekeres vector field of an arbitrary foliation $\mathcal{R}^+_f$. Assume that the GVL class of that foliation  is trivial. Then there exists a $\tilde{\varphi}$-invariant $2$-form $\omega$ on $\Real_+^3$ such that
\begin{equation}\label{Eq4.2}
d\omega = dx_2\wedge dx_1\wedge dx_0.
\end{equation}
Consider the one-parameter group of diffeomorphisms $\varphi_t$ generated by the vector field $V$. Fix some $\xi \in \mathbb{R}$ and define a new $2$-form $\omega'$ in the following way:
\begin{equation*}\label{Eq4.3}
\omega' = \int_\xi^{\xi+1} \tilde{\varphi}_t^*(\omega) dt.
\end{equation*}
\begin{lem}\label{L4.1}
The form $\omega'$ does not depend on $\xi$, it is $\tilde{\varphi}_t$-invariant for every $t \in \mathbb{R}$ and satisfies the equation $d\omega' = dx_2\wedge dx_1\wedge dx_0$.
\end{lem}
{\bf Proof}. The first statement of the lemma is the consequence of the following equality:
\begin{equation*}\label{Eq4.4}
\frac{d}{d\xi} \omega' = \tilde{\varphi}_{\xi+1}^*(\omega) - \tilde{\varphi}_{\xi}^*(\omega)	 = \tilde{\varphi}_{\xi}^*(\tilde{\varphi}^*(\omega)) - \tilde{\varphi}_{\xi}^*(\omega) = \tilde{\varphi}_{\xi}^*(\omega) - \tilde{\varphi}_{\xi}^*(\omega) = 0.
\end{equation*}
Then
\begin{equation*}\label{Eq4.5}
\tilde{\varphi}_t^*(\omega') = \int_\xi^{\xi+1} \tilde{\varphi}_t^*(\tilde{\varphi}_u^*(\omega)) du = \int_\xi^{\xi+1} \tilde{\varphi}_{t+u}^*(\omega) du = \int_{\xi + t}^{\xi + t + 1} \tilde{\varphi}_{v}^*(\omega) dv = \omega'.
\end{equation*}
Finally,
\begin{align*}\label{Eq4.6}
d\omega' = \int_\xi^{\xi+1} d\varphi_t^*(\omega) dt = \int_\xi^{\xi+1} \varphi_t^*(d\omega) dt &= \int_\xi^{\xi+1} \varphi_t^*(dx_2\wedge dx_1\wedge dx_0) dt  \\ =
\int_\xi^{\xi+1} dx_2\wedge dx_1\wedge dx_0\,\,\, dt &= dx_2\wedge dx_1\wedge dx_0.
\end{align*}
This proves the lemma. \qed\vskip1ex

Since $V$ is smooth on $\Real^3$, the new form $\omega'$ is defined and smooth on $\Real_+^3$. Now we may assume that $\omega$ is $\tilde{\varphi}_t$-invariant for every $t \in \mathbb{R}$. Let
\begin{equation*}\label{Eq4.7}
\omega = A dx_2\wedge dx_1 + B dx_0\wedge dx_2 + C dx_1\wedge dx_0,	
\end{equation*}
where $A, B, C: \mathbb{R}_+^3 \rightarrow \mathbb{R}$ are  $C^\infty$-functions of variables $x_0, x_1, x_2$. Consider the vector field $W = (A, B, C)$. The condition (\ref{Eq4.2}) implies
\begin{equation}\label{Eq4.8}
\mbox{div } W = \frac{\partial A}{\partial x_0} + \frac{\partial B}{\partial x_1} + \frac{\partial C}{\partial x_2} = 1.  		
\end{equation}
\begin{lem}\label{L4.2}
	The vector field $U$  of the one-parameter group of diffeomorphisms $\tilde{\varphi}_t$ on $\mathbb{R}^3$ is given by
	\begin{equation*}\label{Eq4.9}
	U = \big(V, V', -x_2 V' + V''\big).
	\end{equation*} 		
\end{lem}

{\bf Proof}. By the definition of $V$, we have
\begin{equation}\label{Eq4.10}
\frac{d}{dt}\big|_{t=0} \varphi_t(x) = V(x).
\end{equation}
Differentiating (\ref{Eq4.10}) by $x$, we obtain
\begin{equation*}\label{Eq4.11}
\frac{d}{dt}\big|_{t=0} \varphi_t'(x) = V'(x),\quad   \frac{d}{dt}\big|_{t=0} \varphi_t''(x) = V''(x).
\end{equation*}
Using \eqref{coordchange}, we get
\begin{multline*}\label{Eq4.12}
U(x_0, x_1, x_2) = \frac{d}{dt}\big|_{t=0} \tilde{\varphi}_t	 (x_0, x_1, x_2)  = \\
\frac{d}{dt}\big|_{t=0} \left( \varphi_t(x_0), x_1 + \ln \varphi_t'(x_0), \frac{x_2}{\varphi_t'(x_0)} + \frac{\varphi_t''(x_0)}{\varphi_t'(x_0)^2} \right)   \\ =
\left( V(x_0), \frac{V'(x_0)}{\varphi_0'(x_0)}, -\frac{x_2 V'(x_0)}{\varphi_t'(x_0)^2} + \frac{V''(x_0)}{\varphi_0'(x_0)^2} -2 \frac{\varphi_0''(x_0) V'(x_0)}{\varphi_0'(x_0)^3}\right). 		
\end{multline*}
It holds $\varphi_0(x) = x$, therefore, $\varphi_0'(x) = 1$ and $\varphi_0''(x) = 0$. This implies
\begin{equation*}\label{Eq4.13}
U(x_0, x_1, x_2) = \big( V(x_0), V'(x_0), -x_2 V'(x_0) + V''(x_0)\big).	
\end{equation*}
\qed\vskip1ex

The $\tilde{\varphi}_t$-invariance of $\omega$ is equivalent to $L_U \omega = 0$, which can be rewritten as
\begin{equation}\label{Eq4.14}
[U, W] = 0.	
\end{equation}

\begin{lem}\label{L4.3}
	Equation (\ref{Eq4.14}) is equivalent to the following system of PDEs:
	\begin{subequations}\label{Eq4.15}
		\begin{align}
		V \frac{\partial A}{\partial x_0} + V' \frac{\partial A}{\partial x_1} + \left( -x_2 V' + V''\right) \frac{\partial A}{\partial x_2}	&= V' A,
		\\
		V \frac{\partial B}{\partial x_0} + V' \frac{\partial B}{\partial x_1} + \left( -x_2 V' + V''\right) \frac{\partial B}{\partial x_2}	&= V'' A,
		\\
		V \frac{\partial C}{\partial x_0} + V' \frac{\partial C}{\partial x_1} + \left( -x_2 V' + V''\right) \frac{\partial C}{\partial x_2}	 &= \left( -x_2 V'' + V'''\right) A - V' C.
		\end{align}
	\end{subequations}
\end{lem}

\noindent{\bf The  proof} may be obtained by the direct computations. \qed\vskip1ex

Up to now we have proved that the  triviality of the GVL class is equivalent to the existence of a common solution $(A, B, C)$ of PDEs (\ref{Eq4.15}) and (\ref{Eq4.8}).
From now on we assume that the foliation is defined by the vector field \eqref{Eq4.1} and we consider the coordinate $x_0$ sufficiently close to $0\in\Real$.
 Substituting  $V$ to (\ref{Eq4.15}), we  get
\begin{subequations}\label{Eq4.16}
	\begin{equation}
	x_0^{2\alpha+2} \frac{\partial A}{\partial x_0} + \alpha x_0^{\alpha+1} \frac{\partial A}{\partial x_1} + \left( - \alpha x_2 x_0^{\alpha+1} -\alpha(\alpha+1) x_0^\alpha + \alpha^2\right) \frac{\partial A}{\partial x_2}	= \alpha x_0^{\alpha+1} A,\end{equation}
		\begin{equation}
	x_0^{2\alpha+2} \frac{\partial B}{\partial x_0} + \alpha x_0^{\alpha+1} \frac{\partial B}{\partial x_1} + \left( - \alpha x_2 x_0^{\alpha+1} -\alpha(\alpha+1) x_0^\alpha + \alpha^2\right) \frac{\partial B}{\partial x_2}
	= \left(  \alpha^2 - \alpha(\alpha+1) x_0^\alpha \right) A,
	\end{equation}
	
	\begin{multline}
	x_0^{3\alpha+3} \frac{\partial C}{\partial x_0} + \alpha x_0^{2\alpha+2} \frac{\partial C}{\partial x_1} + \left( - \alpha x_2 x_0^{\alpha+1} -\alpha(\alpha+1) x_0^\alpha + \alpha^2\right) x_0^{\alpha+1} \frac{\partial C}{\partial x_2}	\\
	= \left( \alpha^3 - 3\alpha^2(\alpha+1)x_0^\alpha + \alpha(\alpha+1)(\alpha+2)x_0^{2\alpha}\right.
	\left. - \alpha^2x_2x_0^{\alpha+1} + \alpha(\alpha+1)x_2x_0^{2\alpha+1} \right) A -
	\alpha x_0^{2\alpha+2} C.
	\end{multline}	
	\end{subequations}\vskip1ex

\noindent Letting $\alpha = n + \beta$ for $n \in \mathbb{Z}$ and $0 \leq \beta <1$, we obtain
\begin{subequations}\label{Eq4.17}
	\begin{equation}
	x_0^{2 n + 2 + 2\beta} \frac{\partial A}{\partial x_0} + \alpha x_0^{n + 1 + \beta} \frac{\partial A}{\partial x_1} + \left( - \alpha x_2 x_0^{n + 1 + \beta} -\alpha(\alpha+1) x_0^{n + \beta} + \alpha^2\right) \frac{\partial A}{\partial x_2}	
	= \alpha x_0^{n + 1 + \beta} A,
	\end{equation}

	\begin{multline}
	x_0^{2 n + 2 + 2\beta} \frac{\partial B}{\partial x_0} + \alpha x_0^{n + 1 + \beta} \frac{\partial B}{\partial x_1} + \left( - \alpha x_2 x_0^{n + 1 + \beta} -\alpha(\alpha+1) x_0^{n + \beta} + \alpha^2\right) \frac{\partial B}{\partial x_2}
\\	= \left(  \alpha^2 - \alpha(\alpha+1) x_0^{n + \beta} \right) A,
	\end{multline}
	
	\begin{multline}
	x_0^{3 n + 3 + 3\beta} \frac{\partial C}{\partial x_0} + \alpha x_0^{2n + 2 + 2\beta} \frac{\partial C}{\partial x_1}
	+ \left( - \alpha x_2 x_0^{2n + 2 + 2\beta} -\alpha(\alpha+1) x_0^{2n + 1 + 2\beta} + \alpha^2 x_0^{n+1+\beta}\right) \frac{\partial C}{\partial x_2}	\\
	= \left( \alpha^3 - 3\alpha^2(\alpha+1)x_0^{n + \beta} + \alpha(\alpha+1)(\alpha+2)x_0^{2n + 2\beta} - \alpha^2 x_2 x_0^{n + 1 + \beta} \right. \\
	\left. + \alpha(\alpha+1)x_2x_0^{2n + 1 + 2\beta} \right) A  - \alpha x_0^{2 n + 2 + 2\beta} C.
	\end{multline}
\end{subequations}
Substituting $x_0=0$ to (\ref{Eq4.17}c), we see that $A(0, x_1, x_2) = 0$. Then it holds  $$A(x_0, x_1, x_2) = x_0 A_1(x_0,x_1, x_2)$$ for some smooth function $A_1$. After dividing the equation (\ref{Eq4.17}c) by $x_0$ we get
\begin{multline}\label{Eq4.18}
x_0^{3 n + 2 + 3\beta} \frac{\partial C}{\partial x_0} + \alpha x_0^{2n + 1 + 2\beta} \frac{\partial C}{\partial x_1}
+ \left( - \alpha x_2 x_0^{2n + 1 + 2\beta} -\alpha(\alpha+1) x_0^{2n + 2\beta} + \alpha^2 x_0^{n + \beta}\right) \frac{\partial C}{\partial x_2}	\\
= \left( \alpha^3 - 3\alpha^2(\alpha+1)x_0^{n + \beta} + \alpha(\alpha+1)(\alpha+2)x_0^{2n + 2\beta} - \alpha^2 x_2 x_0^{n + 1 + \beta} \right. \\
\left. + \alpha(\alpha+1)x_2x_0^{2n + 1 + 2\beta} \right) A_1 - \alpha x_0^{2 n + 1 + 2\beta} C.
\end{multline}
We see that $A_1(0,x_1, x_2) = 0$ and we may assume $A_1 = x_0 A_2$. After doing  $n$ similar steps, we  obtain $$A = x_0^{n+2} \tilde{A}$$ for some  smooth function $\tilde{A}$. The equation (\ref{Eq4.17}c) takes the form
\begin{multline}\label{Eq4.19}
x_0^{2 n + 2 + 3\beta} \frac{\partial C}{\partial x_0} + \alpha x_0^{n + 1 + 2\beta} \frac{\partial C}{\partial x_1}
+ \left( - \alpha x_2 x_0^{n + 1 + 2\beta} -\alpha(\alpha+1) x_0^{n + 2\beta} + \alpha^2 x_0^{\beta}\right) \frac{\partial C}{\partial x_2}	\\
= \left( \alpha^3 - 3\alpha^2(\alpha+1)x_0^{n + \beta} + \alpha(\alpha+1)(\alpha+2)x_0^{2n + 2\beta} - \alpha^2 x_2 x_0^{n + 1 + \beta} \right. \\
\left. + \alpha(\alpha+1)x_2x_0^{2n + 1 + 2\beta} \right) x_0 \tilde{A} - \alpha x_0^{n + 1 + 2\beta} C.
\end{multline}
The condition $\alpha \notin \mathbb{N}$ implies $\beta >0$ and we obtain
\begin{equation}\label{Eq4.20}
\frac{\partial C}{\partial x_2}	(0, x_1, x_2) = 0.	
\end{equation}
Note also that we have already obtained
\begin{equation}\label{Eq4.21}
\frac{\partial A}{\partial x_0}	(0, x_1, x_2) = A_1 (0, x_1, x_2) = 0.	
\end{equation}
Consider now the equation (\ref{Eq4.17}b). Substituting $A=x_0^{n+2} \tilde{A}$, we get
\begin{multline}\label{Eq4.22}
x_0^{2 n + 2 + 2\beta} \frac{\partial B}{\partial x_0} + \alpha x_0^{n + 1 + \beta} \frac{\partial B}{\partial x_1} + \left( - \alpha x_2 x_0^{n + 1 + \beta} -\alpha(\alpha+1) x_0^{n + \beta} + \alpha^2\right) \frac{\partial B}{\partial x_2} \\
= \left(  \alpha^2 x_0^{n+2} - \alpha(\alpha+1) x_0^{2n + 2 + \beta} \right) \tilde{A}.
\end{multline}
This implies that $$\frac{\partial B}{\partial x_2} = x_0 B_1$$ for a smooth function $B_1$. Using similar arguments more $n+1$ times, we obtain ${\partial B}/{\partial x_2} = x_0^{n+2} \tilde{B}$ for some smooth $\tilde B$ and the equation (\ref{Eq4.17}b) takes the form:
\begin{equation}\label{Eq4.23}
\begin{aligned}
x_0^{n + 1 + 2\beta} \frac{\partial B}{\partial x_0} + \alpha x_0^{\beta} \frac{\partial B}{\partial x_1} + \left( - \alpha x_2 x_0^{n + 2 + \beta} -\alpha(\alpha+1) x_0^{n + 1 + \beta} + \alpha^2 x_0\right) \tilde{B} \\
= \left(  \alpha^2 x_0 - \alpha(\alpha+1) x_0^{n + 1 + \beta} \right) \tilde{A},
\end{aligned}
\end{equation}
which implies
\begin{equation}\label{Eq4.24}
\frac{\partial B}{\partial x_1}	(0, x_1, x_2) = 0.	
\end{equation}
Equations (\ref{Eq4.20}), (\ref{Eq4.21}) and (\ref{Eq4.24}) are in contradiction to (\ref{Eq4.8}). This proves the theorem. \qed\vskip.5ex

\begin{cor}
	 If a foliation $\F$ on a  3-dimensional manifold $M$ contains the Reeb component  ${\mathcal R}^+_\alpha$ with $\alpha \notin \mathbb{N}$, then the GVL class of $\F$ is non-trivial.
\end{cor}

Now let the Szekeres field $\tilde{V}$ associated to the diffeomorphism $\psi$ be defined by
\begin{equation}\label{Eq4.25}
\tilde{V}(x) = \left\{
\begin{array}{cc}
e^{-\frac{1}{|x|^\alpha}}, & x< 0,\\
0, & x	\geq 0.
\end{array}
\right.	
\end{equation}
Consider the function $g$ generated by the field $\tilde{V}$ and define the Reeb foliation ${\mathcal R}_\alpha = {\mathcal R}_{g,f}$ depending only on the choice of $\alpha>0$. Theorem \ref{Th4.1} implies that $\GVL({\mathcal R}_\alpha)$ is non-trivial for $\alpha \notin \mathbb{Z}$. The next Theorem shows that the same result is true for any positive odd integer $\alpha$.

\begin{theorem}\label{Th4.2}
	If $\alpha$ is a positive odd  integer, then $\GVL({\mathcal R}_\alpha)$ is non-trivial.
\end{theorem}

{\bf Proof}.  Let $\varphi$ and $\psi$ be the  diffeomorphisms as in Section~\ref{secReeb}. Then the Szekeres field of the diffeomorphism $\varphi \circ \psi^{-1}$ (we denote it by $V(x)$ again) is the following:
\begin{equation}\label{Eq4.26}
V(x) = \left\{
\begin{array}{cc}
-e^{-\frac{1}{|x|^\alpha}}, & x \neq 0,\\
0, & x = 0.
\end{array}
\right.	
\end{equation}
Suppose that $\GVL({\mathcal R}_\alpha)$ is trivial. We get that the equations (\ref{Eq4.15}) and (\ref{Eq4.8}) hold for $V(x)$ defined by (\ref{Eq4.26}).
The proof of the current theorem is similar to the prove of Theorem \ref{Th4.1} with the only difference that we have two sets of equations: for positive and negative $x_0$. We again consider $x_0$ sufficiently close to $0\in\Real$. For $x_0>0$, the equalities (\ref{Eq4.16}) are true again. The analog of (\ref{Eq4.16}) for $x_0<0$ is the following:
\begin{subequations}\label{Eq4.27}
	\begin{equation}
	x_0^{2\alpha+2} \frac{\partial A}{\partial x_0} - \alpha x_0^{\alpha+1} \frac{\partial A}{\partial x_1} + \left( \alpha x_2 x_0^{\alpha+1} +\alpha(\alpha+1) x_0^\alpha + \alpha^2\right) \frac{\partial A}{\partial x_2}	= - \alpha x_0^{\alpha+1} A,
	\end{equation}
	\begin{equation}
	x_0^{2\alpha+2} \frac{\partial B}{\partial x_0} - \alpha x_0^{\alpha+1} \frac{\partial B}{\partial x_1} + \left(  \alpha x_2 x_0^{\alpha+1} +\alpha(\alpha+1) x_0^\alpha + \alpha^2\right) \frac{\partial B}{\partial x_2} \\
	= \left(  \alpha^2 + \alpha(\alpha+1) x_0^\alpha \right) A,
	\end{equation}
	\begin{multline}
	x_0^{3\alpha+3} \frac{\partial C}{\partial x_0} - \alpha x_0^{2\alpha+2} \frac{\partial C}{\partial x_1} + \left(  \alpha x_2 x_0^{\alpha+1} +\alpha(\alpha+1) x_0^\alpha + \alpha^2\right) x_0^{\alpha+1} \frac{\partial C}{\partial x_2}	
	 \\ = \left( - \alpha^3 - 3\alpha^2(\alpha+1)x_0^\alpha - \alpha(\alpha+1)(\alpha+2)x_0^{2\alpha}\right.  \\
	\left. - \alpha^2x_2x_0^{\alpha+1} - \alpha(\alpha+1)x_2x_0^{2\alpha+1} \right) A +
	\alpha x_0^{2\alpha+2} C.
	\end{multline}
	\end{subequations}
Now we may assume that $\beta=1$ and $n \geq 0$, $n \in\mathbb{Z}$. Then from the equation (\ref{Eq4.19}) it follows that
\begin{equation*}\label{Eq4.28}
\frac{\partial C}{\partial x_2} (0, x_1, x_2)= \alpha \tilde{A}(0, x_1, x_2).
\end{equation*}
As in the proof of Theorem \ref{Th4.1}, from  (\ref{Eq4.27}) we get
\begin{equation*}\label{Eq4.29}
\frac{\partial C}{\partial x_2} (0, x_1, x_2)= - \alpha \tilde{A}(0, x_1, x_2).	
\end{equation*}
Therefore we obtain
\begin{equation*}\label{Eq4.30}
\frac{\partial C}{\partial x_2} (0, x_1, x_2) = 0.
\end{equation*}
Similarly,  (\ref{Eq4.23}) implies
\begin{equation*}\label{Eq4.31}
\frac{\partial B}{\partial x_1} (0, x_1, x_2) + \alpha \tilde{B} (0, x_1, x_2) = \alpha \tilde{A} (0, x_1, x_2).
\end{equation*}
For $x_0<0$ the same considerations give
\begin{equation*}\label{Eq4.32}
- \frac{\partial B}{\partial x_1} (0, x_1, x_2) + \alpha \tilde{B} (0, x_1, x_2) = \alpha \tilde{A} (0, x_1, x_2).
\end{equation*}
Consequently,
\begin{equation*}\label{Eq4.33}
\frac{\partial B}{\partial x_1} (0, x_1, x_2) = 0
\end{equation*}
and we obtain a contradiction in the same way as in the proof of Theorem \ref{Th4.1}. \qed

\section[]{Reeb foliations with trivial GVL class}

\begin{theorem}\label{Th5.1}
	If  $\alpha \in \mathbb{N}$, then $\GVL({\mathcal R}^+_\alpha)$ is trivial.
\end{theorem}

{\bf Proof}. We will consider $V$ given by \eqref{Eq4.1} with $\alpha \in \mathbb{N}$ and find a solution to (\ref{Eq4.15}) and (\ref{Eq4.8}) on~$\Real^3_+$.
According to Remark \ref{rem}, this will imply the triviality of~$\GVL({\mathcal R}^+_\alpha)$.

Let  $$\mathbb{R}^3_+ = \{ (x_0, x_1, x_2)\, :\, x_i \in \mathbb{R}, x_0\geq 0\},$$ $$\Gamma = \{ (0, x_1, x_2) \, :\, x_1, x_2 \in \mathbb{R} \} \subset \mathbb{R}^3_+.$$ For any Szekeres field $V$ define the functions $u_1, u_2$ as follows:
\begin{subequations}\label{Eq5.1}
	\begin{align*}
	u_1 &= u_1(x_0,x_1, x_2) = x_1 - \ln |V|, \\
	u_2 &= u_2(x_0,x_1, x_2) =- x_2 V + V'.
	\end{align*}
\end{subequations}
If $V(x_0)\neq0$ for all $x_0>0$, then the functions $u_1, u_2$ define the smooth map
\begin{equation*}
w: \mathbb{R}^3_+ \backslash \Gamma \rightarrow \mathbb{R}^2.
\end{equation*}

To proceed further we need the following lemma.
\begin{lem}\label{L5.1}
	Let $V$ be any Szekeres vector field with $V(x_0)\neq0$ for all~$x_0>0$. Then every pair of smooth functions $\Phi, \Psi: U \rightarrow \mathbb{R}$ defined on an open subset $U \subset \mathbb{R}^2$ gives the following solution to (\ref{Eq4.15}), (\ref{Eq4.8}) smooth on $w^{-1}(U) \subset \mathbb{R}^3_+\backslash \Gamma$:
	\begin{subequations}\label{Eq5.2}
		\begin{align}
		A(x_0, x_1, x_2) &=- V \Phi(u_1, u_2), \\
		B(x_0, x_1, x_2) &= u_1 - \frac{\partial}{\partial u_2}	\Psi(u_1, u_2) - V'\Phi(u_1, u_2), \\
		C(x_0,x_1, x_2) &= -\frac{1}{V}\frac{\partial}{\partial u_1} \Psi(u_1, u_2) + \left(  x_2 V' - V'' \right) \Phi(u_1,u_2).
		\end{align}	
	\end{subequations}	
\end{lem}

\noindent{\bf The proof} of the lemma is by direct substitution of (\ref{Eq5.2}) to (\ref{Eq4.15}) and (\ref{Eq4.8}). \qed\vskip1ex

Now we return  to $V$ defined by (\ref{Eq4.1}). Then
\begin{subequations}\label{Eq5.3}
	\begin{align}
	u_1 &= x_1 + \frac{1}{x_0^\alpha}, \\
	u_2 &= \left(x_2 - \frac{\alpha}{x_0^{\alpha+1}}\right) e^{-\frac{1}{x_0^\alpha}}.
	\end{align}
\end{subequations}

In order to conclude the proof of Theorem~\ref{Th5.1}, we should  find
an example of functions $\Phi$ and~$\Psi$ smooth in the whole~$\Real^2$ such that the corresponding functions $A$, $B$, and $C$ given by~\eqref{Eq5.2} admit a smooth extension from~$\Real^3_+\setminus \Gamma$ to~$\Real^3_+$. To this end, we define
	\begin{subequations}\label{EQ_Phi-Psi}
		\begin{align}
        \Pi &= \{ (u_1, u_2) \, :\, u_1>0, u_2<0 \} \subset \mathbb{R}^2,\notag \\[.75ex]
		\Phi(u_1, u_2) &= -\frac{\big(\alpha+1\big)\,\nu(-2u_1u_2e^{u_1}-1)}{\alpha u_1 u_2} & \text{for~$(u_1,u_2)\in\Pi$},\\
        \Phi(u_1, u_2) &= 0 & \text{for~$(u_1,u_2)\in\Real^2\setminus\Pi$},\\
		\Psi(u_1, u_2) &= u_2 \left(\frac{\alpha + 1}{\alpha} \ln u_1 + 1 - \ln |u_2|\right)\nu(-2u_1u_2e^{u_1}-1) & \text{for~$(u_1,u_2)\in\Pi$},  \\
		\Phi(u_1, u_2) &= 0 & \text{for~$(u_1,u_2)\in\Real^2\setminus\Pi$},
		\end{align}
	\end{subequations}
where $\nu:\Real\to[0,1]$ is a smooth function such that $\nu(t)=0$ for all~$t<0$ and $\nu(t)=1$ for all~${t>1}$.

It is elementary to check that $\Phi$ and $\Psi$ defined above are smooth in~$\Real^2$.
We claim that for any $x_1,x_2\in\Real$ and any $x_0\in\big(0,1/a(x_1,x_2)\big)$, where $$a(x_1,x_2):=\max\Big\{1,|x_1|+\frac{|x_2|+|x_1x_2|+e^{-x_1}}\alpha\Big\},$$ we have $-2u_1u_2e^{u_1}-1\,>\,1$ and hence
\begin{subequations}\label{Eq5.6}
\begin{align}
	A(x_0, x_1, x_2) 	&= \frac{(\alpha+1)x_0^{2\alpha+1}}{\alpha \left( 1 + x_0^\alpha x_1\right) \left( \alpha - x_0^{\alpha+1} x_2 \right)},\\
	B(x_0, x_1, x_2)
	&= x_1 - \ln \frac{\alpha - x_0^{\alpha+1} x_2}{\left(1 + x_0^\alpha x_1 \right)^{1 + \frac{1}{\alpha}}} + \frac{(\alpha + 1) x_0^{\alpha}}{ \left( 1 + x_0^\alpha x_1\right) \left( \alpha - x_0^{\alpha+1} x_2 \right)},\\
	C(x_0, x_1, x_2)
	&= \frac{(\alpha + 1) (\alpha x_0^\alpha x_2  - x_0^{2\alpha+1}x_2^2-\alpha(\alpha+1)x_0^{\alpha-1})}{\alpha \left( 1 + x_0^\alpha x_1\right) \left( \alpha - x_0^{\alpha+1} x_2 \right)},
	\end{align}
\end{subequations}
which can be checked by direct computation using~\eqref{Eq5.2} and~\eqref{Eq5.3}. In particular, this means that $A$, $B$, and $C$ extend smoothly to~$\Real_+^3$ as desired.

It remains to prove our claim. Indeed, fix some $x_1,x_2\in\Real$ and suppose $1/x_0>a(x_1,x_2)$. Then
$$
\frac{\alpha}{x_0^{2\alpha+1}}\ge\frac{\alpha}{x_0^{\alpha+2}}>\frac{\alpha|x_1|+|x_2|+|x_1x_2|+e^{-x_1}}{x_0^{\alpha+1}}
\ge-\frac{\alpha x_1}{x_0^{\alpha+1}}+\frac{x_2}{x_0^{\alpha}}+x_1x_2+e^{-x_1},
$$
which, due to~\eqref{Eq5.3}, implies that $-u_1u_2e^{u_1}>1$ as desired. The proof is now complete. \qed \vskip2ex

\begin{theorem}\label{Th5.2}
	Let  $\alpha \in \mathbb{N}$ be  even. Then $\GVL({\mathcal R}_\alpha)$ is trivial.
\end{theorem}
{\bf Proof.}
It is sufficient to show that for the Szekeres vector field $V$  defined by (\ref{Eq4.26}),  system~\eqref{Eq4.8},\,\eqref{Eq4.15} has a solution~$(A,B,C)$ defined and smooth in the whole~$\Real^3$. The proof of Theorem~\ref{Th5.1} shows that there exists a
smooth solution~$(A,B,C)$ defined on~$\Real_+^3\setminus\Gamma$. Extend the functions $A$, $B$, and~$C$ to~$\Real^3\setminus\Real^3_+$ by setting
\begin{align*}
A(x_0,x_1,x_2)&=-A(-x_0,x_1,-x_2),\\
B(x_0,x_1,x_2)&=\hphantom{-}B(-x_0,x_1,-x_2),\\
C(x_0,x_1,x_2)&=-C(-x_0,x_1,-x_2),\\
\end{align*}
for all $x_0<0$ and all $x_1,x_2\in\Real$.

Using the symmetry of the vector field~$V$, it is easy to see that the extended functions $A$, $B$, and $C$ satisfy system~\eqref{Eq4.8},\,\eqref{Eq4.15} also in~$\Real^3\setminus\Real^3_+$. Moreover, since $\alpha$ is even, equalities~\eqref{Eq5.6} valid for all~$x_1,x_2\in\Real$ and all $x_0\in\big(0,1/a(x_1,x_2)\big)$  extend without any changes to all $(x_0,x_1,x_2)\in\Real^3\setminus\Real^3_+$ with $0>x_0>-1/a(x_1,x_2)$. The right-hand sides of equalities~\eqref{Eq5.6} are smooth at every point of~$\Gamma$. Therefore, the solution $(A, B, C)$ extends smoothly to the whole~$\Real^3$. \qed

\section{Reeb components on surfaces}
The only compact oriented surface that admits a  codimension-one foliation is the torus $T^2$. Any codimension-one foliation on $T^2$ may be constructed in the following way~\cite{Lawson}. Consider an orientation preserving diffeomorphism $h$ of the circle $S^1$.  Consider the  action of the group $\mathbb{Z}$  on $\Real\times S^1$ generated by $(t,\theta)\mapsto(t+1,h(\theta))$. The foliation $\{\Real\times\{\theta\}\} _{\theta\in S^1}$ of $\Real\times S^1$ projects to a foliation on $T^2$. The foliation may be modified by introducing a Reeb component at each fixed point of $h$.

To define a Reeb component, one fixes an even or odd smooth function $f:(-1,1)\to\Real$ strictly monotone on $(-1,0)$ and $(0,1)$ and such that $$\lim_{t\to-1}f(t)=\pm\infty,\quad \lim_{t\to 1}f(t)=\infty,\quad \frac{1}{f'} \text{ is smooth on } [-1,1],$$ considers the foliation of $[-1,1]\times\Real$ by the leaves $x=\pm 1$ and $y=f(x)+c$, $c\in\Real$, and then factorizes this foliation to obtain a foliation on the ring $[-1,1]\times S^1$.

Let us define the structure of the $\D_1$-space on the leaf space of a Reeb component. In a way similar to the one from Section~\ref{secReeb},
consider two transversals $[0,+\infty)\to [-1,1]\times S^1$ (starting from two different circles). We will obtain two holonomy generators $\varphi_1$ and $\varphi_2$ as the local diffeomorphisms of
$[0,+\infty)$ corresponding to the two circles. Since the transversals intersect the same leaves, there appears an orientation preserving diffeomorphism $m$ of $(0,+\infty)$. Thus the leaf space is $$[0,+\infty)/<\varphi_1>\,\, \bigsqcup\,\,  [0,+\infty)/<\varphi_2>$$ factorized by $m$.
In particular, if $f$ is even or odd, then $m=\id_{(0,+\infty)}$, $\varphi_1=\varphi_2=\varphi$. Let the function $f$ be associated to a number $\alpha$ as in the previous sections. We denote the corresponding Reeb component by $\mathfrak{R}_\alpha$. From the above we immediately obtain.

\begin{cor}
	The GVL class of the foliation ${\mathfrak R}_\alpha$ is trivial if and only if $\alpha \in \mathbb{N}$. In particular, if a foliation on a surface contains a Reeb component  	
	${\mathfrak R}_\alpha$ with $\alpha \not\in \mathbb{N}$, then its GVL class is non-trivial.
\end{cor}

\begin{cor} If $\alpha \in \mathbb{N}$ and $\beta \notin \mathbb{N}$, then the
	Reeb components ${\mathfrak R}_\alpha$ and ${\mathfrak R}_\beta$ are not diffeomorphic. 	
\end{cor}

\vskip0.1cm

{\bf Acknowledgements.} The authors are thankful to Steven Hurder
for useful email communications. The work was supported by grant no.
18-00496S of the Czech Science Foundation. Ya.\,V.\,Bazaikin was partially supported by the Program of Fundamental Scientific Research of the SB RAS No. I.1.2., Project No. 0314-2019-0006.

\end{document}